\documentclass[11pt,reqno]{amsart}
\usepackage{amsmath,amssymb,mathrsfs}
\usepackage{fullpage}
\usepackage{xcolor}
\usepackage[colorlinks=true,allcolors=blue,pagebackref]{hyperref}
\hypersetup{pdftitle={The monodromy conjecture for remarkable poles},
 pdfauthor={Nero Budur, Quan Shi, Huaiqing Zuo},pdfsubject={Research draft},
 pdfkeywords={Bernstein-Sato polynomial, log canonical modification, monodromy}}
\numberwithin{equation}{section}
\newtheorem{theorem}{Theorem}[section]
\newtheorem{lemma}[theorem]{Lemma}
\theoremstyle{definition}
\newtheorem{remark}[theorem]{Remark}
\newcommand{\C}{\mathbb C}
\newcommand{\Q}{\mathbb Q}
\newcommand{\Z}{\mathbb Z}
\newcommand{\OO}{\mathcal O}
\newcommand{\Dscr}{\mathscr D}
\DeclareMathOperator{\ord}{ord}
\DeclareMathOperator{\Supp}{Supp}
\DeclareMathOperator{\Res}{Res}
\DeclareMathOperator{\Fix}{Fix}
\DeclareMathOperator{\lct}{lct}
\title[Monodromy for remarkable poles]{The monodromy conjecture for remarkable poles}
\author{Nero Budur}
\address{Department of Mathematics, KU Leuven, Celestijnenlaan 200B, 3001 Leuven, Belgium;  YMSC, Tsinghua University, 100084 Beijing, China;  BCAM, Mazarredo 14, 48009 Bilbao, Spain.}
\email{nero.budur@kuleuven.be}

\author{Quan Shi}
\address{Department of Mathematical Sciences, Tsinghua University, Beijing, 100084, P. R. China.}
\email{shiq24@mails.tsinghua.edu.cn / thusq20@gmail.com}

\author{Huaiqing Zuo}
\address{Department of Mathematical Sciences, Tsinghua University, Beijing, 100084, P. R. China.}
\email{hqzuo@mail.tsinghua.edu.cn}

\begin{document}
\begin{abstract} We prove the strong monodromy conjecture for the remarkable numbers of a regular function $f$ on a smooth complex variety. We show more generally that the log canonical modification of $f$ computes the remarkable numbers and gives roots of the Bernstein-Sato polynomial of $f$.
\end{abstract}
\maketitle
\tableofcontents

\section{Introduction}\label{sec:intro}

Let $X$ be a smooth irreducible complex variety and let $f\in\OO_X(X)$
be nonzero with nonempty divisor $D=\operatorname{div}(f)$. The strong
monodromy conjecture predicts that if $s_0$ is a pole of the motivic
zeta function of $f$, then $b_f(s_0)=0$, where $b_f$ is the
Bernstein--Sato polynomial. The monodromy conjecture predicts that
$\exp(2\pi i s_0)$ is a monodromy eigenvalue on the nearby cycles of
$f$. Budur, de Lorenzo Poza, Shi, and Zuo singled out a set of
positive rational numbers $\alpha$, called remarkable numbers, and
proved that the negative of each remarkable number is a motivic pole
\cite[Theorem~1.2(2)]{BDPSZ}. Our main result proves both predictions
for these poles, namely \cite[Conjecture~1.3]{BDPSZ}.

Take a log resolution $\mu:V\to X$ that is an isomorphism off $D$, and
write $\mu^*D=\sum_{i\in J}N_iE_i$ and
$K_{V/X}=\sum_{i\in J}(\nu_i-1)E_i$. For each nonempty $I\subset J$
with $\bigcap_{i\in I}E_i\ne\varnothing$, put
$N_I=\gcd_{i\in I}N_i$ and $\alpha_I=\min_{i\in I}\nu_i/N_i$.
The $d$-log canonical threshold, for an integer $d\ge1$, is
\begin{equation}\label{eq:d-definition}
 {}^{(d)}\lct(f)=\min_{I:\,N_I\mid d}\alpha_I,\qquad \min\varnothing=+\infty .
\end{equation}
For fixed $d$ the value is independent of the resolution, and its finite values as $d$ varies are precisely the remarkable numbers \cite{BDPSZ}. Taking $d$ divisible by all $N_i$ gives the usual log canonical threshold.

\begin{theorem}\label{thm:main}
If $\alpha$ is a remarkable number of $f$, then $b_f(-\alpha)=0$.
In particular, $\exp(-2\pi i\alpha)$ is an eigenvalue of the monodromy on the nearby cycles of
$f$. Equivalently, it is an eigenvalue on $H^j(F_y,\C)$ for some
$y\in D$ and some $j$, where $F_y$ is the Milnor fiber at $y$.
\end{theorem}

Here $b_f(s)$ is the monic polynomial of least degree for which
$b_f(s)f^s\in\Dscr_X[s]f^{s+1}$ locally on $X$, with $\Dscr_X$ the
sheaf of algebraic differential operators. The polynomial obtained using local analytic differential operators
 at $x$ is denoted by $b_{f,x}(s)$. It is well-known that
$
 b_f(s)=\operatorname{lcm}_{x\in X(\C)} b_{f,x}(s).
$
%Indeed, in \'etale coordinates, prescribing the order and $s$-degree of a Bernstein operator for a given polynomial makes its identity a finite linear system over $\OO_{X,x}$, after clearing a power of $f$ and comparing coefficients of $s$. For $b_{f,x}$ this system has a solution in the analytic local ring, hence also in $\OO_{X,x}$ by faithful flatness. The resulting identity holds on a Zariskineighborhood of $x$. Finitely many such neighborhoods cover $X$, so the least common multiple of their polynomials gives an identity locally on $X$. Conversely, every algebraic identity analytifies at each point.

\subsection{Divisors on the log canonical model}
Let \(\pi:S\to X\) be the log canonical modification of \((X,D_{\mathrm{red}})\), that is, the log canonical model over \(X\) of \((Y,(h^*D)_{\mathrm{red}})\) for any log resolution \(h:Y\to X\) that is an isomorphism over \(X\setminus D\).
Thus $\pi$ is projective and birational,
$B_S=(\pi^*D)_{\mathrm{red}}$ is the sum of the strict boundary
and the exceptional prime divisors, $(S,B_S)$ is log canonical, and
$K_S+B_S$ is ample over $X$. We call a prime component of $B_S$ a
\emph{retained divisor}. For a prime divisor $G$ over $X$, let $\ord_G$
be its normalized order of vanishing and $\nu_G=A_X(G)$ its log
discrepancy. If $N_G=\ord_G(f)>0$, its slope is $\nu_G/N_G$.

\begin{theorem}[Log canonical model witnesses]\label{thm:witness}
Every remarkable number of $f$ is the slope of a retained divisor.
More precisely, suppose $X$ is quasi-projective, take a small
projective $\Q$-factorialization $r:T\to S$, and let the $B_j$ be the
strict transforms of the components of $B_S$. By smallness these are
the primes of $((\pi r)^*D)_{\mathrm{red}}$. There is an integer $L_0>0$
such that for every positive multiple $L$ of $L_0$ the following holds.
Let $\rho_T:T_L\to T$ be the normalized cyclic cover $t^L=f\circ\pi\circ r$. Let
$\mu_L$,  the group of \(L\)-th roots of unity, act by scalar multiplication on $t$ and choose a generator
$\zeta\in\mu_L$. For all $d\ge1$,
\begin{equation}\label{eq:fixed-locus}
 {}^{(d)}\lct(f)=
 \min_{\;j:\,B_j\cap\rho_T(\Fix(\zeta^d;T_L))\ne\varnothing}
          \frac{\nu_{B_j}}{N_{B_j}} .
\end{equation}
\end{theorem}

A single point of $B_j$ with a $\zeta^d$-fixed lift to $T_L$ suffices
for $B_j$ to enter the minimum in \eqref{eq:fixed-locus}.

\begin{theorem}[Bernstein--Sato roots from the log canonical model]
\label{thm:retained}
If $G$ is a retained divisor, then
\[
       0<\frac{\nu_G}{N_G}\le1,\qquad
       b_f\!\left(-\frac{\nu_G}{N_G}\right)=0 .
\]
If $\pi(G)=\{x\}$, the same root occurs in $b_{f,x}(s)$.
\end{theorem}

Together, Theorems~\ref{thm:witness} and~\ref{thm:retained} give the
root in Theorem~\ref{thm:main}. The Malgrange--Kashiwara correspondence
then gives the monodromy eigenvalue.

\medskip
\noindent\textbf{Overview.} The proof of Theorem \ref{thm:witness} is given in Section~\ref{sec:witness}. 
On a log resolution, a point on precisely the components indexed by $I$
has a $\zeta^d$-fixed lift to a sufficiently divisible cyclic cover
exactly when $N_I\mid d$. Mapping this lift to $T_L$ and comparing
discrepancies gives one inequality in \eqref{eq:fixed-locus}.
For the converse, Section~\ref{sec:witness} starts from a fixed lift
over $B_j$. Every prime extracted by a crepant divisorial log
terminal modification of $(T,B_j)$ has the slope of $B_j$. Fixed-point lifting on
the normalized cover and on a minimal log canonical center then gives
a point satisfying the gcd condition on a divisor of that slope.

For Theorem~\ref{thm:retained}, the strategy pushes a classical approach, with a new twist. First suppose $G$ has center a point
$x$. Put $\alpha=\nu_G/N_G$ and let $E$ be its strict transform on a
log resolution $\mu:Y\to X$ dominating $S$. For a local volume form
$\omega$, the residue $\rho=\Res_E((f\circ\mu)^{-\alpha}\mu^*\omega)$
has rank-one unitary local coefficients. 
Logarithmic Akizuki--Nakano vanishing gives
$[\rho]\ne0$ in the hypercohomology of a suitable logarithmic de Rham
complex on $E$ in Lemma~\ref{lem:first-lattice}.
The crucial ingredient here that allows one to apply this vanishing is Lemma \ref{lem:signed-ample} providing a suitable ample divisor; it rests on the existence of log resolutions that are isomorphisms over the simple normal crossings locus of the pair and uses the log canonical model assumption. 

With the class $[\rho]$ in hand, the rest of the approach is classical, and consists of using the integration pairing to obtain Bernstein-Sato roots. 
By logarithmic duality and Serre vanishing, Lemma~\ref{lem:test}
constructs a meromorphic function $h$ near $x$ whose pullback restricts
to a rational function $h_E$ on $E$ satisfying
$\langle[\rho],\overline{[h_E\rho]}\rangle\ne0$. This pairing can be expressed in terms of the residue at $-\alpha$ an archimedean zeta integral weighted by $\overline h$, by Lemma \ref{lem:pairing}, which is more general and we think is of independent interest, and Lemma \ref{lem:zeta-residue}. The Bernstein identity introduces a factor
of $f$, making the shifted integral holomorphic at $-\alpha$.
Lemma~\ref{lem:point-root} thus gives $b_{f,x}(-\alpha)=0$.
For a positive-dimensional center, Lemma~\ref{lem:slice-root} transfers
the root from general transverse slices to $b_f$.

\medskip
\noindent\textbf{Acknowledgements.} We thank M. \v Cmrlec for useful discussions; he declined co-authorship of this paper. The authors were supported by grants KUL  Methusalem METH/21/03, FWO  G0B3123N, BJNSF  1252009, NSFC  125B2004, NSFC  12671056. LLMs available with a ChatGPT Pro subscription were used.

\medskip
\noindent\textbf{AI use.} AI was essential. We had attempted many times unsuccessfully with and without  AI to prove Theorem \ref{thm:main}. Meantime, from discussions with M. \v Cmrlec it emerged without AI that his work on a different topic implies that remarkable numbers are computed on dlt modifications. We prompted AI to assume this. It gave first a proof conditional on remarkable numbers being computed on log canonical modifications, then answered positively our prompt to prove Theorem \ref{thm:witness}. We verified all AI-generated content.

\section{The cyclic fixed locus and the log canonical model}
\label{sec:witness}

We prove Theorem~\ref{thm:witness} computing ${}^{(d)}\lct(f)$ from  divisors on log canonical models.

For a normal pair $(W,\Gamma)$ with $K_W+\Gamma$ $\Q$-Cartier and a
prime $G$ on a birational model $q:V\to W$, its log discrepancy is
$A_{W,\Gamma}(G)=1+\ord_G(K_V-q^*(K_W+\Gamma))$, with compatible
canonical divisors. The pair is log canonical (lc) if these numbers
are nonnegative, and Kawamata log terminal (klt) if they are positive.
The centers of divisors of log discrepancy zero are the lc centers.
A divisorial log terminal (dlt) pair has simple normal crossings at
the generic points of its lc centers.

A point on precisely the components $E_i$, $i\in I\ne\varnothing$,
of a log resolution is \emph{$d$-admissible} if $N_I\mid d$.
Each nonempty simple normal crossings intersection has a dense open
stratum, so \eqref{eq:d-definition} minimizes $\alpha_I$ over
$d$-admissible closed points.

For a normal birational model $W\to X$, let $\rho_W:W_L\to W$ be its
normalization in $\C(X)[t]/(t^L-f)$, a product of fields. It is the
normalization of the cover $t^L=f$, whose algebra is locally free and
embeds in its generic fiber. Normalization is finite for complex
varieties. The action of $\mu_L$ on $t$ and proper birational maps
$W'\to W$ lift to these normalizations because elements integral over $\OO_W$
pull back to the integrally closed local rings on the corresponding
components of $W'_L$. The induced maps $W'_L\to W_L$ are equivariant,
birational on corresponding components and compatible with composition
by uniqueness. They are isomorphisms over the isomorphism locus of
$W'\to W$, and proper because $W'_L\to W$ is proper and $W_L\to W$
is separated.

\begin{lemma}\label{lem:cyclic-fiber}
Suppose that \(V\) is  smooth as in  \S\ref{sec:intro} and \(\operatorname{div}_V(f)\) has simple normal crossing
support. Let \(v\) be a closed point lying on precisely the components
\(E_i\), \(i\in I\ne\varnothing\). The geometric fiber of \(V_L\) at
\(v\), as a set, is a transitive \(\mu_L\)-set of cardinality
\[
h=\gcd\bigl(L,N_{E_i}:i\in I\bigr).
\]
In particular, if \(N_{E_i}\mid L\) for all \(i\in I\), the point
\(v\) is \(d\)-admissible if and only if it has a lift fixed by
\(\zeta^d\), for any generator \(\zeta\in\mu_L\).
\end{lemma}

\begin{proof}
In étale crossing coordinates write $f=u\prod_{i\in I}x_i^{N_{E_i}}$.
An étale extension adjoining an $L$-th root of the unit $u$ absorbs
it into $t$. These base changes preserve the geometric fiber and its
action, since the base-changed normalization is finite, birational and
normal. The equation now factors into the primitive binomials
\[
 t^{L/h}-\xi\prod_{i\in I}x_i^{N_{E_i}/h},\qquad \xi\in\mu_h .
\]
After rescaling $t$ in one factor, its monomials $x_i,t$ generate the
lattice $M=\Z^I+\Z(N_{E_i}/L)_{i\in I}$ and span the positive
orthant. The normalized semigroup is therefore $M\cap\mathbb R_{\ge0}^I$. The other coordinates are smooth factors. Every nonconstant monomial
in the saturation has a positive power in the ideal $(x_i:i\in I)$,
so the normalized branch has one point above the crossing. Multiplication
by $\zeta$ permutes the branches by the primitive $h$-th root
$\zeta^{L/h}$, giving the assertion.
\end{proof}

\begin{lemma}\label{lem:fixed-lift}
Let \(a\) be a closed point of a complex variety \(Z\) that has
rational singularities near \(a\). Let \(g\) be an automorphism of
finite order fixing \(a\). For any proper equivariant birational
morphism \(Z'\to Z\), there is a \(g\)-fixed point of \(Z'\) above
\(a\).
\end{lemma}

\begin{proof}
The assertion is immediate in dimension zero. Otherwise restrict to
the intersection of the $g$-translates of a neighborhood of $a$ with
rational singularities. Resolve $Z'$ functorially and principalize the
pullback of the ideal of $a$ with simple normal crossings support
\cite[Theorems~3.35 and~3.36]{Kol}. This gives a smooth proper
birational model $q:\widetilde Z\to Z$ dominating $Z'$, with a
$g$-action and reduced simple normal crossings fiber $P=q^{-1}(a)_{\rm red}$.
Since $Rq_*\OO_{\widetilde Z}=\OO_Z$, \cite[Lemma~2.5]{MS} gives
$H^0(P,\OO_P)=\C$ and $H^k(P,\OO_P)=0$ for $k>0$. The alternating
trace of $g$ is therefore one.

If $P$ had no fixed point, the alternating resolution of $\OO_P$ by the
structure sheaves of the intersections of its components would give
trace zero. Permuted summands contribute zero, and on each preserved
intersection holomorphic Lefschetz for its Dolbeault complex gives zero
\cite[Theorem~2]{AB}, regardless of the sign induced by permuting its
indices. Thus $g$ fixes a point of $P$, whose image in $Z'$ lies above $a$.
\end{proof}

\begin{proof}[Proof of Theorem~\ref{thm:witness}]
Take a finite affine cover of $X$. The relative log canonical models
exist on its members and restrict to relative log canonical models on
intersections \cite[Theorem~1.1 and Proposition~2.1]{OX}. Uniqueness
over the base gives the cocycle condition, so the models glue to
$\pi:S\to X$. This morphism is proper locally on $X$, hence proper.
A common Cartier multiple of $K_S+B_S$ is relatively ample on the
cover. A sufficiently high common power $M$ of its line bundle gives
a closed immersion $S\hookrightarrow\mathbb P_X(\pi_*M)$, as this
holds on each member by relative very ampleness. Restricting one log
resolution shows that ${}^{(d)}\lct(f)$ is the minimum of the thresholds on
the cover, with value $+\infty$ on opens disjoint from $D$. A prime
realizing a finite minimum on one chart closes to a prime on $S$ with
the same valuation. We may therefore work over quasi-projective $X$.

By uniqueness, $\pi$ is an isomorphism over the smooth complement of
$D$. For sufficiently small rational $\epsilon>0$, the boundary
$B_S-\epsilon\pi^*D$ is effective and the resulting pair is klt. Its log
discrepancies increase by $\epsilon\ord_G(f)>0$ for divisors centered
in $B_S$, and the complement is smooth. By \cite[Theorem~10.4]{Fuj}
there is a small projective $\Q$-factorialization. Take any such
$r:T\to S$, set $h=\pi r$, and write $B_T=(h^*D)_{\rm red}=\sum_jB_j$.
Smallness identifies the boundary primes and gives the crepant identity
$K_T+B_T=r^*(K_S+B_S)$ in codimension one, hence as divisors. In
particular $(T,B_T)$ is lc and $(T,B_T-\epsilon h^*D)$ is klt.
Removing its effective $\Q$-Cartier boundary gives that $(T,0)$ is klt.
A small projective birational map to a normal $\Q$-factorial variety
is an isomorphism. Indeed, a multiple of a relatively ample divisor equals the
pullback of its pushforward, so no curve can be contracted. This applies
to $r$ over $S\setminus B_S$.

For a divisor $E$ over $T$, write $a_j(E)=\ord_E(B_j)$. Pulling back
$h^*D$ and $K_T+B_T-h^*K_X=\sum_j\nu_{B_j}B_j$ gives
\begin{equation}\label{eq:witness-discrepancy}
N_E=\sum_j a_j(E)N_{B_j},\qquad
\nu_E=A_{T,B_T}(E)+\sum_j a_j(E)\nu_{B_j}.
\end{equation}
For each $B=B_j$, take a $\Q$-factorial crepant dlt modification
$p:(Y,C)\to(T,B)$ by the same theorem of Fujino. Its exceptional
primes $G$ have coefficient one in $C$ and satisfy
\[
0\le A_{T,B_T}(G)=A_{T,B}(G)-\ord_G(B_T-B)=-\ord_G(B_T-B)\le0.
\]
Since $(T,0)$ is klt, each such $G$ is centered in $B$. Thus $H=p^*(B_T-B)$ is the strict transform of $B_T-B$, and
$(Y,C+H)$ is crepant over $(T,B_T)$. The description of the
exceptional primes and \eqref{eq:witness-discrepancy} give
\begin{equation}\label{eq:witness-equal-slope}
C=(p^*B)_{\rm red},\qquad \frac{\nu_G}{N_G}=\frac{\nu_B}{N_B}
\quad\text{for every component }G\text{ of }C,
\end{equation}
since an exceptional $G$ has $A_{T,B_T}(G)=0$ and $a_i(G)=0$ for
$B_i\ne B$. Over $T\setminus B$ the map $p$ is small and hence an
isomorphism by the preceding $\Q$-factoriality argument.

Choose a log resolution $V\to T$ of $(T,B_T)$ that is an isomorphism
off $B_T$. For each $(Y,C+H)$, apply \cite[Theorem~10.45(2)]{KolSing}
to choose a log resolution $V_B\to Y$ that is an isomorphism over the
simple normal crossings locus of $(Y,C+H)$. The composites to $X$ are
log resolutions of $f$, since all the models are isomorphic off $D$.
Choose $L_0$ divisible by all multiplicities on $V$ and
the finitely many $V_B$. For $L_0\mid L$ and $d\ge1$, choose a generator
$\zeta\in\mu_L$ and put $g=\zeta^d$ and $F_d=\rho_T(\Fix(g;T_L))$.

If $v\in V$ is $d$-admissible, Lemma~\ref{lem:cyclic-fiber} and
equivariance send a fixed lift of $v$ to a fixed lift of its image
$z\in T$. For a component $E$ through $v$, the indices with
$a_j(E)>0$ form a nonempty set by \eqref{eq:witness-discrepancy}.
Each such $B_j$ contains the center of $E$ and therefore contains
$z$, the image of $v\in E$. The fixed lift of $z$ gives $z\in F_d$,
so $z\in B_j\cap F_d$.
The same identity and $A_{T,B_T}(E)\ge0$ give
\[
\min_{j:B_j\cap F_d\ne\varnothing}\frac{\nu_{B_j}}{N_{B_j}}
\le\min_{j:a_j(E)>0}\frac{\nu_{B_j}}{N_{B_j}}\le\frac{\nu_E}{N_E}.
\]
Taking the minimum on $V$ gives the first inequality in \eqref{eq:fixed-locus}.

Conversely, choose $B$ meeting $F_d$ and a fixed closed point $z_L$
above $B$. For boundary primes $P$ on $T$ or on the corresponding $Y$,
write $\bar P$ for the entire reduced inverse image in the corresponding normalized cyclic cover, and extend the bar
additively. In the strict henselization at $P$ the unit in $f=u x^{N_P}$
has an $L$-th root. As in Lemma~\ref{lem:cyclic-fiber}, there are
$\gcd(L,N_P)$ primitive binomial factors, each of ramification index
$e_P=L/\gcd(L,N_P)$, unchanged by strict henselization. Hence
$\bar P=e_P^{-1}\rho^*P$ is $\Q$-Cartier. The covers are étale off
the boundary, so tame ramification \cite[(2.41.4)]{KolSing} gives
\[
 K_{T_L}+\overline{B_T}=\rho_T^*(K_T+B_T),\qquad
 K_{Y_L}+\bar C+\bar H=\rho_Y^*(K_Y+C+H).
\]
The finite-cover discrepancy formula \cite[Proposition~5.20(4) and its proof]{KM} gives
log canonicity of both pairs. Subtracting $\epsilon\operatorname{div}_{T_L}(f)$
from $\overline{B_T}$ gives the log crepant pullback of the klt pair
$(T,B_T-\epsilon h^*D)$ and an effective boundary after further
decreasing $\epsilon>0$. Thus $T_L$ has rational singularities
\cite[Theorem~5.22]{KM}. The component containing $z_L$ is preserved
by $g$, since $T_L$ is normal. Lemma~\ref{lem:fixed-lift} applied to
the corresponding component of $Y_L\to T_L$ gives a fixed point $y_L$
over $z_L$, lying in $\bar C$ by \eqref{eq:witness-equal-slope}.

The pair $(Y_L,\bar C)$ is lc after removing the effective $\Q$-Cartier
divisor $\bar H$. A coefficient-one prime passes through $y_L$, so
there is a unique minimal lc center $A$ through it
\cite[Theorem~9.1(1),(4)]{Fuj}, necessarily $g$-invariant. Remove the
finitely many lc centers strictly contained in $A$, which all miss
$y_L$, and choose an affine neighborhood of $y_L$ in the remaining
open. There $A$ is a minimal lc center, since every lc center of the
restriction closes to one of the original pair. Local subadjunction
\cite[Theorem~7.2]{FG} gives rational singularities of $A$ near $y_L$.

Since $(Y_L,\bar C+\bar H)$ is lc, $A$ is not contained in
$\Supp\bar H$. Over $Y\setminus\Supp H$ the cover is log crepant
for $\bar C$ and $C$. A divisorial valuation of log discrepancy zero
centered at $A$ restricts, after normalization, to a divisorial
valuation centered at $\rho_Y(A)$. Its log discrepancy downstairs
is zero by the finite-cover discrepancy formula. Thus $\rho_Y(A)$
is an lc center of $(Y,C)$. Its generic point lies on a component
$G$ of $C$ and outside $\Supp H$. There $(Y,C+H)$ has simple normal
crossings because $(Y,C)$ is dlt. Hence $V_B\to Y$ is an isomorphism
over a neighborhood of that point, and $(V_B)_L\to Y_L$ is an isomorphism
over the generic point of $A$. The strict transform $A'$ of $A$ is
$g$-invariant and proper birational over $A$. The image of its generic point in $V_B$ lies
on the strict transform of $G$, so the whole image lies there. Applying
Lemma~\ref{lem:fixed-lift} to $A'\to A$ at $y_L$ yields a fixed
point whose image in $V_B$ is $d$-admissible by
Lemma~\ref{lem:cyclic-fiber}. The component through this image coming
from $G$ has slope $\nu_B/N_B$ by \eqref{eq:witness-equal-slope}.
Hence ${}^{(d)}\lct(f)\le\nu_B/N_B$. This proves \eqref{eq:fixed-locus},
including the empty case, and realizes each finite threshold on $T$,
hence on $S$ by smallness.
\end{proof}

\section{Regularized integration for complementary extensions}
\label{sec:regularized-pairing}

The normal and tangential integrations in Lemma~\ref{lem:zeta-residue}
use the next two lemmas. The first concerns archimedean local zeta
distributions in normal crossing coordinates. The one-variable complex
powers are classical \cite[Appendix~B, \S B1.3]{GS}. An early
independent-parameter calculation for real coordinate monomials is
\cite[Lemma~2]{BG}.

\begin{lemma}[Archimedean local zeta integral in normal crossing coordinates]\label{lem:local-mellin}
Let $z_1,\ldots,z_e$ be coordinates on a polydisc and $0\le r\le e$.
For $1\le i\le r$, fix $a_i,b_i\in\mathbb R$ with $a_i-b_i\in\Z$.
Write $dA(z)=\tfrac i2\,dz\wedge d\bar z$ and interpret
$z^a\bar z^{\,b}$ as $|z|^{a+b}(z/|z|)^{a-b}$ on the punctured disc.
For a compactly supported smooth function $\phi$, put
\[
 I(\mathbf v;\phi)=\int \phi(z)
   \prod_{i=1}^r |z_i|^{2v_i}z_i^{a_i}\bar z_i^{\,b_i}
                  \frac{dA(z_i)}{|z_i|^2}
   \prod_{j=r+1}^e dA(z_j),\qquad \Re v_i\gg0 .
\]
As a distribution in $z$, this extends meromorphically in the independent
parameters $\mathbf v=(v_1,\ldots,v_r)$. Its only possible poles lie on
\begin{equation}\label{eq:angular-poles}
 v_i=-\max(a_i,b_i)-k,\qquad k\in\Z_{\ge0},
\end{equation}
with order at most one along each hyperplane. If $a_i=b_i=0$, then
\[
 \Res_{v_i=0} I(\mathbf v;\phi)
 =\pi\, I(\mathbf v_{\widehat i};\phi|_{z_i=0}),
\]
where the right side omits the $i$-th factor.
These formulas also hold for $\phi$ depending holomorphically on
$\mathbf v$ and auxiliary parameters, with locally common compact support
and locally uniform bounds on every coordinate derivative. In the residue
formula such a coefficient is restricted to $v_i=0$ as well.
\end{lemma}

\begin{proof}
For one coordinate, take $f(z)=z$ in Igusa's Theorem~5.3.2 and the
remark following it \cite[\S5.3, p.~72]{Ig}. With respect to $dA$, our density
has exponents $s_1=v+a-1$ and $s_2=v+b-1$ in $z^{s_1}\bar z^{s_2}$.
Thus Igusa's angular character has exponent $s_1-s_2=a-b$.
Since $b_z(s)=s+1$, division of this distribution by either
$\Gamma(v+a)$ or $\Gamma(v+b)$ gives an entire family. Its poles
therefore lie in the intersection of the two gamma pole sets, giving
\eqref{eq:angular-poles} with order at most one. Igusa applies this
calculation on a resolution in \cite[proof of Theorem~5.4.1]{Ig}.

To keep the normal parameters independent, apply the same integration
by parts in each coordinate. Put
$K_i=|z_i|^{2v_i-2}z_i^{a_i}\bar z_i^{\,b_i}$ and $m_i=\max(a_i,b_i)$.
Take $(\eta_i,D_i)=(z_i,\partial_{z_i})$ if $a_i\ge b_i$ and
$(\bar z_i,\partial_{\bar z_i})$ otherwise. For $k\ge1$ write
$(x)_k=x(x+1)\cdots(x+k-1)$. Off $z_i=0$ one has
$D_i^k(\eta_i^k K_i)=(v_i+m_i)_k K_i$. For all $\Re v_i$ sufficiently
large the boundary terms at $z_i=0$ vanish, giving
\[
 \Bigl(\prod_{i=1}^r(v_i+m_i)_k\Bigr) I(\mathbf v;\phi)
 =(-1)^{rk}\int \Bigl(\prod_{i=1}^r\eta_i^k K_i\Bigr)
    (D_1^k\cdots D_r^k\phi)\prod_{j=1}^e dA(z_j).
\]
The right side is jointly holomorphic for $2\Re v_i+a_i+b_i+k>0$
for every $i$. On compact subsets of this domain, the absolute values
of the kernels have a common integrable majorant on the support of
$\phi$, and the integral is bounded by a smooth seminorm of order $rk$.
This also proves joint holomorphy for the stated families of test
functions. Increasing $k$ covers parameter space. Each denominator
has distinct roots in its own variable, proving the pole and order
assertions for distributions.

For $a_i=b_i=0$, start with the other real parts large and subtract
from $\phi$ its restriction at $z_i=0$ extended constantly in $z_i$
and multiplied by a smooth cutoff equal to one on $|z_i|<R$.
The difference is $O(|z_i|)$ there, so its integral is holomorphic
near $v_i=0$. The cutoff annulus also contributes a holomorphic term.
With our area form the integral on the disc is
\[
 \int_{|z_i|<R}|z_i|^{2v_i-2}\,dA(z_i)=\frac{\pi R^{2v_i}}{v_i}.
\]
Taking the residue restricts the coefficient to $v_i=0$. Meromorphic
continuation in the remaining parameters gives the stated formula.
\end{proof}

Let $E$ be a smooth projective complex variety of dimension $m$,
$B=\sum_i B_i$ a reduced simple normal crossings divisor, and
$U=E\setminus B$. Fix a rank-one local system $\mathcal L$ on $U$
with a flat Hermitian metric $h$. A logarithmic lattice $M$ for $\mathcal L$
is a line bundle on $E$ with an integrable logarithmic connection whose
horizontal local system on $U$ is identified with $\mathcal L$. Write
$C(M)=(\Omega_E^\bullet(\log B)\otimes M,\nabla)$ for its analytic
de Rham complex in degrees $0,\ldots,m$.

\begin{lemma}[Integration and logarithmic duality]\label{lem:pairing}
Let $M_a,M_b$ be logarithmic lattices for $\mathcal L$ with residues
$a_i,b_i\in\mathbb R$ along $B_i$. Thus $a_i-b_i\in\Z$, since the two local monodromies coincide.
Assume that whenever $a_i$ is integral, exactly one of $a_i,b_i$ is
positive. Let $\xi,\eta$ be closed forms of total degree $m$ in the
smooth Dolbeault resolutions of $C(M_a),C(M_b)$, respectively. For
each $i$, give $\OO_E(B_i)$ a smooth Hermitian metric, let $s_i$ be its
canonical section, that is, the image of $1$ under the natural inclusion
$\OO_E\hookrightarrow\OO_E(B_i)$, so that $\operatorname{div}(s_i)=B_i$,
and take $c_i>0$. Then
\[
 \operatorname{Reg}\!\int_E\xi\wedge\overline\eta
 :=\left.\left(\int_U\prod_i\|s_i\|^{2c_iw}\,
              \xi\wedge\overline\eta\right)\right|_{w=0}
\]
is defined by meromorphic continuation from $\Re w\gg0$, holomorphic
at zero. The wedge uses the flat metric on $\mathcal L$.
This value is independent of the metrics on $\OO_E(B_i)$ and the $c_i$.
On cohomology it gives the perfect pairing
\[
 \mathbb H^m(E,C(M_a))\otimes_{\C}
       \overline{\mathbb H^m(E,C(M_b))}\longrightarrow\C
\]
given by logarithmic Verdier duality.
\end{lemma}

\begin{proof}
\emph{The dual logarithmic extension.}
Let $N=M_a^\vee(-B)$ with the connection induced by the dual connection
and ordinary differentiation on $\OO_E(-B)$. Its local system on $U$
is $\mathcal L^\vee$. Its residue along $B_i$ is $1-a_i$, since a
local generator of $\OO_E(-B)$ has a simple zero there.
Esnault--Viehweg's logarithmic duality
\cite[Appendix~A, Proposition~A.2]{EV86} gives
\[
 C(N)\simeq\mathbb D_E C(M_a)[-2m],
 \quad \text{ where }\ \mathbb D_E F=R\mathcal Hom_{\C}(F,\C_E[2m]).
\]
The shift uses the complex orientation of $E$.

Let $P$ be the union of the $B_i$ for which $b_i$ is a positive
integer. The hypothesis says precisely that $1-a_i$ is a positive
integer on the same components. Write
$U\xrightarrow{i}E\setminus P\xrightarrow{j}E$ for the open inclusions.
By Esnault--Viehweg's logarithmic comparison
\cite[\S1.4 and Lemma~1.5]{EV89}, both complexes are obtained by
derived direct image across $B\setminus P$ and extension by zero across
$P$:
\begin{equation}\label{eq:uniqueness_CMb}
	\overline{C(M_b)}\simeq j_!Ri_*\overline{\mathcal L},
	\qquad C(N)\simeq j_!Ri_*\mathcal L^\vee .
\end{equation}
Here the bar conjugates the sheaf of complex vector spaces. The flat
metric identifies $\overline{\mathcal L}$ with $\mathcal L^\vee$ on
$U$. Since $Ri_*$ and $j_!$ are fully faithful for open immersions,
this identification extends uniquely to an isomorphism
$\overline{C(M_b)}\simeq C(N)$. It remains to identify its trace with
the integral in the statement.

\emph{The weighted integration current.}
We construct a family of currents by integrating the Hermitian wedge
product against a weight that ensures convergence near $B$.
For $\ell=a,b$, let $\mathscr A_\ell^\bullet$ be the total smooth
Dolbeault resolution of $C(M_\ell)$, with differential $D_\ell$.
On $U$ this differential is the exterior covariant derivative.
Let $W\subset E^{\rm an}$ be open and take
$u\in\mathscr A_a^p(W)$ and $v\in\mathscr A_b^q(W)$. On $W\cap U$ the
flat metric $h$ on the common flat line bundle, linear in its first
argument, defines a complex-valued
$(p+q)$-form. For local flat sections $s,t$ of $\mathcal L$, set
\[
 \beta(\alpha\otimes s,\overline{\gamma\otimes t})
 =h(s,t)\,\alpha\wedge\overline\gamma.
\]
In a common unit flat frame $e$ (so $s=t=e$),
this is $\alpha\wedge\overline\gamma$.
Put $k=p+q$ and $w_{\mathbf t}=\prod_i\|s_i\|^{2t_i}$ for independent
complex parameters $\mathbf t=(t_i)_i$. For $k\leq 2m$ and
$\Re t_i\gg0$, define a degree-$k$ current on compactly supported
smooth $(2m-k)$-forms $\varphi$ on $W$ by
\[
 T_{\mathbf t}(u,\bar v)(\varphi)
 =\int_{W\cap U}w_{\mathbf t}\,\beta(u,\bar v)\wedge\varphi .
\]
We use zero currents in degrees above $2m$.

\emph{Continuation to the origin.}
We show that $T_{\mathbf t}(u,\bar v)$ continues holomorphically to a
neighborhood of $\mathbf t=0$ and that its value there commutes with
restriction to open subsets.
Choose coordinates on a polydisc in which the boundary components
meeting it are $B_i=(z_i=0)$ for $1\le i\le r$. Let $e_a,e_b$ be local
holomorphic frames of $M_a,M_b$. On the complement of the boundary,
choose a common multivalued horizontal frame $e$ of unit norm for the
flat Hermitian metric $h$. Then
\[
 e_a=f_a\prod_{i=1}^r z_i^{a_i}e,\qquad
 e_b=f_b\prod_{i=1}^r z_i^{b_i}e,
\]
where $f_a,f_b$ are holomorphic units on the polydisc. Indeed,
subtracting the residue terms from either connection form gives a
closed holomorphic one-form on the polydisc, hence an exact form.
Thus the coefficient $h(e_a,e_b)$ in $\beta(u,\bar v)$ is
$f_a\overline{f_b}\prod_{i=1}^r z_i^{a_i}\bar z_i^{\,b_i}$,
which is single-valued since $a_i-b_i\in\Z$.

The top-degree form $w_{\mathbf t}\beta(u,\bar v)\wedge\varphi$ on the complement
of the boundary has the form
\[
 \left(\prod_{i=1}^r |z_i|^{2t_i}z_i^{a_i}\bar z_i^{\,b_i}\right)
 F_{\mathbf t}(z)
 \left(\bigwedge_{i=1}^r\frac{dA(z_i)}{|z_i|^2}\right)
 \wedge
 \left(\bigwedge_{j=r+1}^m dA(z_j)\right),
\]
with $F_{\mathbf t}$ smooth across the boundary. Indeed, each nonzero
term contains at most one $dz_i/z_i$ and one $d\bar z_i/\bar z_i$.
A missing logarithmic denominator contributes a factor $z_i$ or
$\bar z_i$ to $F_{\mathbf t}$. The dependence on $\mathbf t$ is
holomorphic since $\|s_i\|^2=|z_i|^2e^{\psi_i}$ for $i\le r$, with
$\psi_i$ smooth, and the remaining sections are nowhere zero here.

A partition of unity on the support of $\varphi$ now reduces the
integral to Lemma~\ref{lem:local-mellin}. Its seminorm estimates give
continuity in $\varphi$ after continuation. By
\eqref{eq:angular-poles}, a pole through $\mathbf t=0$ would require
$a_i,b_i\in\mathbb Z_{\leq0}$ for some $i$, which the hypothesis
excludes. Thus $T(u,\bar v):=T_{\mathbf 0}(u,\bar v)$ is a current.
It commutes with restriction since extension by zero of a test form
preserves the defining integral and hence its continuation.

\emph{Independence of the regularization.}
We show that the value at zero is independent of the boundary metrics
and is recovered along every line $t_i=c_iw$ with $c_i>0$.
Multiplying the squared boundary norms by $e^{h_i}$, for smooth real
functions $h_i$, replaces $\varphi$ by $e^{\sum_i t_i h_i}\varphi$.
Lemma~\ref{lem:local-mellin} applies to this holomorphic family of test
forms, whose value at zero is $\varphi$. Hence $T$ is independent of the
boundary metrics. For $\mathbf c=(c_i)_i$ with every $c_i>0$, the
restriction $\mathbf t=w\mathbf c$ gives the one-parameter integral
for $\Re w\gg0$ and has the same value $T(u,\bar v)$ at $w=0$ by
joint holomorphy. For the dependence of the finite part on
the regulator in the presence of poles, see
\cite[Theorem~5.4(2) and Corollary~5.5]{FK}.

\emph{Compatibility with the differentials.}
We verify the graded product rule using the convention
$(dS)(\varphi)=(-1)^{k+1}S(d\varphi)$ for a degree-$k$ current.
If $2\Re t_i+a_i+b_i>3$ for every $i$, the coefficients of
$w_{\mathbf t}\beta(u,\bar v)$ extend by zero as $C^1$ functions
across $B$. Indeed, the logarithmic differentials contribute at most
two powers of $|z_i|^{-1}$ in each boundary coordinate, and one
derivative costs at most one further power. By flatness of $h$ and the identity
$dw_{\mathbf t}=w_{\mathbf t}\sum_i t_i\,d\log\|s_i\|^2$ on $U$,
Stokes' formula on $W$ gives
\begin{equation}\label{eq:Stokes}
 dT_{\mathbf t}(u,\bar v)
 =T_{\mathbf t}(D_a u,\bar v)
  +(-1)^pT_{\mathbf t}(u,\overline{D_b v})
  +\sum_i t_i S_{i,\mathbf t},
\end{equation}
where $S_{i,\mathbf t}$ is the current defined by integration against
$w_{\mathbf t}\,d\log\|s_i\|^2\wedge\beta(u,\bar v)$. For the
one-parameter identity, see
\cite[equation~(3.1) in the proof of Proposition~3.2(ii)]{FK}.

Near $B_i=(z_i=0)$ we have
$d\log\|s_i\|^2=dz_i/z_i+d\bar z_i/\bar z_i+d\psi_i$.
Since repeated differentials wedge to zero, insertion of this factor
preserves the local form of the top-degree integrand above. By \eqref{eq:angular-poles} in
Lemma~\ref{lem:local-mellin}, the possible poles of $S_{i,\mathbf t}$
lie on $t_j=-\max(a_j,b_j)-k$ with $k\in\Z_{\geq0}$. Since
$a_j-b_j\in\Z$, such a hyperplane contains the origin only if
$a_j,b_j\in\Z_{\leq0}$, which the residue hypothesis excludes.
Thus $S_{i,\mathbf t}$ is holomorphic at $\mathbf t=0$. Evaluating the continued Stokes identity \eqref{eq:Stokes} at $\mathbf t=0$ gives
\begin{equation}\label{eq:product_rule}
	dT(u,\bar v)
	=T(D_a u,\bar v)+(-1)^pT(u,\overline{D_b v}).
\end{equation}

\emph{Identification with the Verdier duality pairing.}
We identify the morphism induced by $T$ with logarithmic duality and
compute its orientation trace on $\xi$ and $\eta$.
Write $\mathscr C_E^\bullet$ for the de Rham complex of currents on
$E^{\rm an}$. The product rule \eqref{eq:product_rule} makes $T$ a morphism of complexes
$\mathscr A_a^\bullet\otimes_{\C}\overline{\mathscr A_b^\bullet}
\to\mathscr C_E^\bullet$. Since currents resolve $\C_E$ by
\cite[Proposition~B.4]{FSY}, tensor--Hom
adjunction with the graded symmetry gives a morphism
\[
 \overline{C(M_b)}\longrightarrow R\mathcal Hom_{\C}(C(M_a),\C_E)
 \simeq\mathbb D_E C(M_a)[-2m].
\]
On $U$ the current is ordinary integration of the wedge, so on
horizontal sections this morphism is the metric identification
$\overline{\mathcal L}\to\mathcal L^\vee$. By \eqref{eq:uniqueness_CMb} and the full faithfulness of $j_!Ri_*$,
this morphism is the unique extension of the metric identification to $E$
and hence is the logarithmic duality isomorphism. For the global closed
degree-$m$ forms $\xi,\eta$ in the statement, the orientation trace
of their product is the top-degree current $T(\xi,\bar\eta)$
evaluated on the constant test function $1$ on compact $E$. By
construction this is $\operatorname{Reg}\!\int_E\xi\wedge\overline\eta$.
By \cite[Proposition~C.4 and Corollary~C.6]{FSY}, composition with this
trace gives the perfect pairing of Verdier duality.
\end{proof}

\begin{remark}
The identification of the two perfect pairings seems new in this generality, but it is not unexpected.
For a rank-one unitary local system on a punctured projective line with
nontrivial monodromy at every puncture, Deligne--Mostow
\cite[\S2.18 and Proposition~2.19]{DM} express the cup product pairing
on classes represented by forms of the first kind or their antiholomorphic
counterparts as the convergent integral of their Hermitian wedge product.
For constant coefficients and the complementary twists by reduced
boundary divisors, the integral converges without regularization.
Brown--Dupont identify it with the relative cohomology pairing
in \cite[Theorem~3.1 and Proposition~3.14]{BD}. For higher-order poles,
they discuss how the cutoff integral can depend on the shape of the
excised neighborhoods in \cite[Remark~3.22]{BD}.
\end{remark}

\section{A residue on a retained divisor}\label{sec:residue-lattices}

For a retained divisor $G_0$ over a closed point $x$, we construct the nonzero residue pairing \eqref{eq:test-pairing}.

For every retained divisor $G$,
\begin{equation}\label{eq:slope-bound}
       0<\nu_G\le\ord_G(D_{\mathrm{red}})\le N_G .
\end{equation}
Indeed, the $\pi$-exceptional divisor
$K_S+B_S-\pi^*(K_X+D_{\mathrm{red}})$ is $\Q$-Cartier and
$\pi$-nef. By the negativity lemma \cite[Lemma~3.39(1)]{KM} its
coefficients are nonpositive. Its coefficient at $G$ is
$\nu_G-\ord_G(D_{\mathrm{red}})$. Positivity of $\nu_G$ follows
from smoothness of $X$.

Shrink to a smooth affine neighborhood of $x$ on which $\omega_X$ has a
nowhere vanishing section $\omega$. Choose compatible canonical divisors,
using $\omega$ on $X$. Take a projective log resolution $p:Y\to S$ of
$(S,B_S)$, write $\mu=\pi p$, and let $E$ be the strict transform of
$G_0$. Write $E_0=E$ and $\sum_{i=0}^r E_i=(\mu^{-1}(f=0))_{\rm red}$.
Set
\[
 \alpha=\nu_0/N_0,\qquad \delta_i=\nu_i-\alpha N_i,\qquad
 m=\dim E,\qquad B_i=E_i|_E\quad(i>0),\qquad B=\sum_{i>0}B_i.
\]
Then $E$ is smooth projective and $B$ is a reduced simple normal
crossings divisor. Each $B_i$ may be empty or disconnected. A coefficient on $B_i$
applies to all of its components.

\begin{lemma}\label{lem:signed-ample}
The resolution $p$ can be chosen to admit a $\mu$-ample rational divisor
$T=\sum_{i=0}^r t_iE_i$ satisfying
\[
 t_0=0,\qquad t_i<\delta_i\ \text{ and }\ t_i\ne0\quad(1\leq i\leq r).
\]
\end{lemma}

\begin{proof}
Choose $p$ by the construction in the proof of
\cite[Theorem~10.45(2)]{KolSing}, omitting blowups of Cartier divisors,
which are isomorphisms. The remaining blowup centers lie over the
complement of the simple normal crossings locus of $(S,B_S)$. This
complement has codimension at least two since $S$ is normal and
$B_S$ is reduced. At each blowup
the negative of the effective tautological divisor is relatively ample. Summing their
pullbacks, with each earlier coefficient sufficiently large relative
to the later ones, gives an effective $p$-exceptional divisor $H$ with
$-H$ ample over $S$ (take $H=0$ for the identity).

The divisor
\[
 Q=K_S+B_S-\pi^*K_X-\alpha\operatorname{div}_S(f)
\]
is $\pi$-ample. Write $p^*Q=\sum_i q_i E_i$. The discrepancy formula gives
$\delta_i-q_i=A_{S,B_S}(E_i)\geq0$. In particular $q_0=\delta_0=0$.
For sufficiently small rational $\eta>0$, the divisor
$H'=-H-\eta\sum_{i>0} E_i$ is $p$-ample. It has coefficient zero on $E_0$
and strictly negative coefficients on all other $E_i$. For sufficiently
small rational $\epsilon>0$, the divisor $T=p^*Q+\epsilon H'$ is
$\mu$-ample. Avoiding finitely many values of $\epsilon$ gives all the
assertions.
\end{proof}

Put $U=E\setminus B$, $\mathcal K=\omega_E(B)$, and write
$\alpha=p_0/q$ in lowest terms. As a rational section of
$\omega_Y(\sum_i E_i)^{\otimes q}$ the tensor
\[
        (f\circ\mu)^{-p_0}(\mu^*\omega)^{\otimes q}
\]
has divisor $\sum_i q\delta_i E_i$, since in a logarithmic volume
frame $\mu^*\omega$ has order $\nu_i$ along $E_i$ and is a unit off
the boundary. As $\delta_0=0$, restriction and logarithmic adjunction
give a nonzero rational section $\tau$ of $\mathcal K^{\otimes q}$.
Its logarithmic coefficient in crossing coordinates is a unit times
$\prod_{i>0}z_i^{q\delta_i}$ over the other components through the
point. Thus $\operatorname{div}_E(\tau)=q\sum_{i>0}\delta_iB_i$.

For a local algebraic frame $\kappa$ of $\mathcal K$, write
$\tau=F_\kappa\kappa^{\otimes q}$ and $e=\kappa^\vee$. On
$\mathcal L_0=\mathcal K^{-1}$ set
\[
 \nabla e=\tfrac1q\,d\log F_\kappa\otimes e,\qquad
 \|e\|^2=|F_\kappa|^{2/q}\quad\text{on }U .
\]
Under $\kappa'=v\kappa$ one has $e'=v^{-1}e$ and
$F_{\kappa'}=v^{-q}F_\kappa$. These formulas give a global integrable
algebraic logarithmic connection with residues $\delta_i$, and a
flat Hermitian metric with local horizontal unit frames
$F_\kappa^{-1/q}e$ on $U$. The section
$\rho=\kappa\otimes e=1$ of $\mathcal K\otimes\mathcal L_0$ has scalar
realizations $F_\kappa^{1/q}\kappa$ in those frames. These are precisely
the local branches of
$\Res_E((f\circ\mu)^{-\alpha}\mu^*\omega)$: the form has a simple
logarithmic pole in the normal coordinate, since $\alpha N_0=\nu_0$,
and the $q$-th tensor power of its residue is $\tau$.

For an integral divisor $G=\sum_{i>0}g_i B_i$, put
$\mathcal L_G=\mathcal L_0(G)$ with its induced logarithmic connection.
Its residues are $\delta_i-g_i$, and on $U$ it has the same
horizontal local system $\mathcal L=\ker(\nabla|_U)$ and flat metric
as $\mathcal L_0$. Its algebraic de Rham complex
$K_G^\bullet$ has terms $K_G^a=\Omega_E^a(\log B)\otimes\mathcal L_G$.
For $\iota:E^{\rm an}\to E_{\rm Zar}$, analytification of forms commutes
with the connection, giving a map $\iota^{-1}K_G^\bullet\to
K_{G,{\rm an}}^\bullet$ of complexes of sheaves of $\C$-vector spaces.
The induced comparison for the filtration by form degree is an
isomorphism on $E_1^{a,b}$ by coherent GAGA for the terms $K_G^a$ on
projective $E$. Boundedness thus identifies
$\mathbb H^m(E,K_G^\bullet)$ with
$\mathbb H^m(E^{\rm an},C(\mathcal L_G^{\rm an}))$, carrying the class of
a global logarithmic top form to the class of its analytification.

For the lattice used below put $c_i=\lfloor\delta_i-t_i\rfloor$ for
$i>0$ and $C=\sum_{i>0}c_iB_i$. Explicitly, in degrees $0\le a\le m$,
\[
 K_C^a=\Omega_E^a(\log B)\otimes\omega_E^{-1}
                    \otimes\OO_E(C-B),\qquad K_C^m=\OO_E(C),
\]
and $K_C^a=0$ in all other degrees. To write the differential, let
$g$ be a local equation for $C$. In the frame $g^{-1}\kappa^\vee$
of $\mathcal L_C$, it sends a logarithmic $a$-form $\theta$ to
\[
 d\theta+\bigl(\tfrac1q\,d\log F_\kappa-d\log g\bigr)\wedge\theta .
\]

\begin{lemma}\label{lem:first-lattice}
The section $\rho$ defines a nonzero class in
$\mathbb H^m(E,K_C^\bullet)$.
\end{lemma}

\begin{proof}
By Lemma~\ref{lem:signed-ample}, $t_i<\delta_i$ for $i>0$, so
$c_i=\lfloor\delta_i-t_i\rfloor\geq0$ and $C$ is effective. Set $L=K_E+B$ and
\[
 \Delta=\sum_{i>0}(\delta_i-t_i-c_i)B_i ,
 \qquad A=L-C-\Delta .
\]
The divisor of $\tau$, a rational section of $\OO_E(qL)$, gives
$L\sim_{\Q}\sum_{i>0}\delta_iB_i$. Since $t_0=0$ by
Lemma~\ref{lem:signed-ample}, $A\sim_{\Q}\sum_{i>0}t_iB_i=T|_E$.
This restriction is ample because $T$ is $\mu$-ample by the same lemma
and $\mu(E)=\{x\}$. By the definition of $c_i$, all coefficients of
$\Delta$ lie in $[0,1)$, so
$\lceil A\rceil=L-C$ and $\lceil A\rceil-A=\Delta\leq B$.

Logarithmic Akizuki--Nakano vanishing
\cite[Theorem~1.3]{AMPW}, applied to $A$ with logarithmic
boundary and integral divisor both equal to $B$, gives
\[
 H^{m-b}\!\left(E,\Omega_E^{m-a}(\log B)\otimes
       \OO_E(\lceil A\rceil-B)\right)=0
 \qquad(a+b<m).
\]
The logarithmic wedge pairing gives
\[
 (K_C^a)^\vee\otimes\omega_E
  \simeq \Omega_E^{m-a}(\log B)\otimes\omega_E(-C)
  =\Omega_E^{m-a}(\log B)\otimes\OO_E(\lceil A\rceil-B).
\]
Serre duality therefore identifies the dual of the vanishing group
with $H^b(E,K_C^a)$. This use of logarithmic vanishing is inspired
by \cite[Lemma~4.4]{Wu}.

In the spectral sequence $E_1^{a,b}=H^b(E,K_C^a)$, all terms of
total degree $m-1$ vanish. For each $\ell\geq1$ the differential into $E_\ell^{m,0}$
has source $E_\ell^{m-\ell,\ell-1}=0$, and the outgoing differential
has target beyond the last column. Hence $E_\infty^{m,0}=E_1^{m,0}$.
The section $1\in H^0(E,\OO_E(C))$ is nonzero and survives to
$\mathbb H^m(E,K_C^\bullet)$, where it represents $\rho$.
\end{proof}

\begin{remark}[Topological description of the logarithmic complex]
Work in the analytic topology. The residue of $\mathcal L_C$ along
$B_i$ is $a_i=\delta_i-\lfloor\delta_i-t_i\rfloor$.
Since $a_i=\lceil t_i\rceil$ when $\delta_i\in\Z$, the union of
components with positive integral residue is
\[
 P=\bigcup_{a_i\in\Z_{>0}}B_i
   =\bigcup_{\delta_i\in\Z,\ t_i>0}B_i .
\]
Put $B_*=\bigcup_{a_i\notin\Z_{>0}}B_i$, $V=E\setminus B_*$
and $A=P\cap V$, so that $U=V\setminus A$. For the open inclusions
\[
 U\xrightarrow{i}E\setminus P\xrightarrow{j}E,
 \qquad U\xrightarrow{\ell}V\xrightarrow{k}E,
\]
logarithmic comparison \cite[\S1.4 and Lemma~1.5]{EV89} gives
\[
 K_{C,{\rm an}}^\bullet
 \simeq j_!Ri_*\mathcal L
 \simeq Rk_*\ell_!\mathcal L .
\]
Thus the boundary condition is extension by zero across $P$ and
derived direct image across $B_*$. Along each component of $A$ the
local monodromy is $\exp(-2\pi\mathrm{i}a_i)=1$. Hence
$\mathcal L$ extends to the local system
$\widetilde{\mathcal L}=\ell_*\mathcal L$ on $V$.
For the closed inclusion $h:A\hookrightarrow V$, the exact sequence
\[
 0\longrightarrow\ell_!\mathcal L\longrightarrow\widetilde{\mathcal L}
 \longrightarrow h_*(\widetilde{\mathcal L}|_A)\longrightarrow0
\]
identifies $H^r(V,\ell_!\mathcal L)$ with relative cohomology.
Consequently, for every $r$,
\[
 \mathbb H^r(E,K_C^\bullet)\simeq H^r(V,A;\widetilde{\mathcal L}).
\]
Under this identification, restriction to $H^r(U,\mathcal L)$ is
induced by the inclusion of pairs $(U,\varnothing)\hookrightarrow(V,A)$.
\end{remark}

\begin{remark}
For plane curves, Loeser detects the residue on a rupture divisor by a
nonzero twisted period on the complement of its boundary points
\cite[proof of Theorem~III.3.1(a)]{Loe88}. Under his nonresonance
and nonvanishing hypotheses, Blanco detects higher coefficients in the
asymptotic expansion of Milnor-fiber periods by analogous twisted periods
\cite[(13.1) and Proposition~13.4]{Bla21}. For Newton-nondegenerate
hypersurfaces, Loeser
proves that the residue on a compact facet divisor has nonzero class in
the ordinary cohomology of its open stratum with local coefficients
when the boundary monodromies are nontrivial
\cite[Theorem~3.7 and proof of Theorem~4.2, \S\S4.3--4.4]{Loe90}.

Here Lemma~\ref{lem:first-lattice} proves $[\rho]\ne0$ in
$\mathbb H^m(E,K_C^\bullet)$ for the logarithmic lattice determined by
$T$, but nonvanishing of its image under restriction to $H^m(U,\mathcal L)$
could in fact fail, in contrast with the above works. However, Lemma~\ref{lem:test} detects
$[\rho]$ by pairing it with a class in the complementary logarithmic
lattice on $E$, and this will be used to provide a root of the Bernstein-Sato polynomial.
\end{remark}

\begin{lemma}\label{lem:test}
For a sufficiently large and divisible integer $k$ there is a rational
function $h=P/f^a$, with $a\geq0$ and $P\in\OO_X(X)$, satisfying
\[
 \ord_{E_0}(h)=0,\qquad \ord_{E_i}(h)\geq-kt_i\quad(i>0).
\]
Write $h_E$ for its restriction at the generic point of $E$ and
$G_k=(kT)|_E$. In the pairing of Lemma~\ref{lem:pairing} for
$\mathcal L_C$ and $\mathcal L_{G_k}$ one has
\begin{equation}\label{eq:test-pairing}
 \bigl\langle[\rho],\overline{[h_E\rho]}\bigr\rangle
 =\operatorname{Reg}\!\int_E\rho\wedge\overline{h_E\rho}\ne0 .
\end{equation} For every $i>0$, the two exponents
$\delta_i$ and $\delta_i+\ord_{E_i}(h)$ cannot both belong to
$\Z_{\leq0}$.
\end{lemma}

\begin{proof}
Choose $k_0>0$ with $k_0T$ integral and $\OO_Y(k_0T)$ very ample
over $X$. Its restriction to $E\subset\mu^{-1}(x)$ is very ample.
Since $t_i\ne0$ for $i>0$ by Lemma~\ref{lem:signed-ample},
take a sufficiently large multiple $k$ of $k_0$ with
$k>|\delta_i/t_i|$ for all $i>0$. The residues of $\mathcal L_C$
and $\mathcal L_{G_k}$ on $B_i$ are
\[
 a_i=\delta_i-\lfloor\delta_i-t_i\rfloor,
 \qquad b_i=\delta_i-kt_i.
\]
For $\delta_i\in\Z$ one has $a_i=\lceil t_i\rceil$. If $t_i>0$, then
$a_i>0$ and $b_i<0$. If $t_i<0$, then $a_i\le0$ and $b_i>0$.
For $\delta_i\notin\Z$, both residues are nonintegral since $kt_i\in\Z$.
Thus the two lattices satisfy Lemma~\ref{lem:pairing}. By the comparison
above and Lemma~\ref{lem:first-lattice}, there exists
$\gamma\in\mathbb H^m(E,K_{G_k}^\bullet)$ with
$\langle[\rho],\overline\gamma\rangle\ne0$.

Apply Serre vanishing \cite[III, Theorem~5.2(b)]{Har} on $E$ to the
finitely many sheaves $\Omega_E^a(\log B)\otimes\mathcal L_0$
and on $Y$, projective over affine $X$, to $\OO_Y(-E)$, twisting
by the corresponding powers of $\OO_Y(k_0T)$. For all sufficiently
large multiples $k$ this gives
\[
 H^b(E,K_{G_k}^a)=0\quad(b>0),\qquad H^1(Y,\OO_Y(kT-E))=0.
\]
By the first vanishing, $\gamma$ is computed by the complex of global
sections. Its degree-$m$ cohomology is the quotient of
$H^0(E,K_{G_k}^m)$ by the image of the preceding differential
(with zero image if $m=0$). Since $K_{G_k}^m=\OO_E(G_k)\rho$, choose
$h_E\in H^0(E,\OO_E(G_k))$ with $[h_E\rho]=\gamma$.
Lemma~\ref{lem:pairing} now gives \eqref{eq:test-pairing} for these
representatives. The second vanishing, applied to the restriction
sequence for $E$, makes $H^0(Y,\OO_Y(kT))\to H^0(E,\OO_E(G_k))$
surjective. A lift $h$ has $\ord_E(h)=0$ because $t_0=0$ and $h_E\ne0$,
and $\ord_{E_i}(h)\ge-kt_i$ since $\operatorname{div}_Y(h)+kT\ge0$.
For $t_i<0$ this gives $\delta_i+\ord_{E_i}(h)\ge\delta_i-kt_i>0$,
while for $t_i>0$ one has $\delta_i>t_i>0$. Finally $\mu$ is an
isomorphism and $T$ vanishes over $X_f$, so $h\in\OO_X(X)[1/f]$.
\end{proof}

\section{The residue of a local zeta integral}\label{sec:mellin}

We prove Theorem~\ref{thm:retained} for a divisor centered at a point
by computing a zeta residue.

Use the notation of Section~\ref{sec:residue-lattices} and take
$h=P/f^a$ from Lemma~\ref{lem:test}. Put $n=\dim X$ and
$c_j=2^{-j}i^{j^2}$, so $c_n\omega\wedge\bar\omega$ is the product
of Euclidean area forms when $\omega=dz_1\wedge\cdots\wedge dz_n$.
Take an analytic coordinate neighborhood $U$ of $x$ small enough for
the local Bernstein identity and choose $\chi\in C_c^\infty(U)$
equal to one near $x$.
For $\Re s\gg0$ define
\begin{equation}\label{eq:mellin-integral}
 Z_h(s)=c_n\int_U |f|^{2s}\overline h\,\chi\,
                 \omega\wedge\overline\omega .
\end{equation}
On $U\setminus D$ one has
$|f|^{2s}\overline h=f^s\bar f^{\,s-a}\overline P$. Thus $Z_h$ is
the archimedean local zeta integral with Igusa's exponents
$(s_1,s_2)=(s,s-a)$ and smooth test density
$c_n\overline P\chi\,\omega\wedge\overline\omega$
\cite[\S5.3, definition following Theorem~5.3.2]{Ig}.

For plane curves, Lichtin \cite[proof of Theorem~2, (2.1)--(2.2)]{Lic}
and Blanco \cite[Proposition~5.4, (5.8), and Remark~5.6]{Bla}
express archimedean zeta residues as integrals along exceptional curves.
Here the choice of $h$ isolates $E$ and realizes the nonzero pairing
\eqref{eq:test-pairing} as its contribution, up to the factor below.

\begin{lemma}[Residue on the retained divisor]\label{lem:zeta-residue}
The meromorphic continuation of $Z_h$ has a simple pole at $s=-\alpha$ with
\begin{equation}\label{eq:mellin-residue}
 \Res_{s=-\alpha}Z_h(s)
  =\frac{\pi}{N_E}\,c_{n-1}\,\operatorname{Reg}\!
     \int_E\rho\wedge\overline{h_E\rho}\ \ne 0.
\end{equation}
\end{lemma}

\begin{proof}
Put $w=s+\alpha$. Cover the compact set $\mu^{-1}(\Supp\chi)$ by
finitely many crossing charts on $Y$ and choose a subordinate
partition of unity. For each boundary prime $G$
through the center of a chart let $z_G$ be its coordinate and
$\ell_G=\ord_G(h)$. Write
\[
 f\circ\mu=u\prod_G z_G^{N_G},\qquad
 \mu^*\omega=u'\prod_G z_G^{\nu_G}\Omega,\qquad
 h\circ\mu=H\prod_G z_G^{\ell_G},
\]
where $u,u'$ are holomorphic units and $\Omega$ is the coordinate
logarithmic volume form. The coefficient $H$ is holomorphic because
$h$ has poles only on the boundary and all boundary orders have been
factored out. The chart contribution to $Z_h(s)$ is
the specialization $v_G=N_Gw$ of the integral in
Lemma~\ref{lem:local-mellin} with factors
\begin{equation}\label{eq:mellin-factors}
 |z_G|^{2v_G}z_G^{\delta_G}\bar z_G^{\delta_G+\ell_G}
 \frac{dA(z_G)}{|z_G|^2},
 \qquad \delta_G=\nu_G-\alpha N_G .
\end{equation}
The smooth coefficient contains $|u|^{2(w-\alpha)}|u'|^2\overline H$
and the cutoffs. It is entire in $w$, which we first keep independent
of the $v_G$. Lemma~\ref{lem:local-mellin} gives meromorphic
continuation on the line $v_G=N_Gw$, since $N_G>0$ for every $G$.

By Lemma~\ref{lem:test}, the exponents in \eqref{eq:mellin-factors}
for $G\ne E$ cannot both be nonpositive integers. Lemma~\ref{lem:local-mellin}
therefore gives holomorphy near $(\mathbf v,w)=0$ on charts disjoint
from $E$. On a chart meeting $E$, one has $\delta_E=\ell_E=0$, so the
only possible pole there is simple along $v_E=0$. Multiplying the chart
integral by $v_E$ gives a holomorphic germ whose restriction to $v_E=0$
is $\pi$ times the tangential integral on $z_E=0$. Since $v_E=N_Ew$
on the specialization, the residue in $s$ is $\pi/N_E$ times the value
of this tangential integral at $(\mathbf v_{\widehat E},w)=0$.

To identify that value, put $t=z_E$ first in the logarithmic volume
form. On a local branch off the other boundary components write
$(f\circ\mu)^{-\alpha}\mu^*\omega=(dt/t)\wedge\beta$. The restriction
$\beta|_E$ is the scalar realization of $\rho$ in a horizontal unit
frame, and $(h\circ\mu)|_E=h_E$ as rational functions since $\ell_E=0$.
The identity $c_n(-1)^{n-1}=(i/2)c_{n-1}$ shows that taking out
$dA(t)/|t|^2$ leaves precisely
$c_{n-1}\rho\wedge\overline{h_E\rho}$ with the flat metric, times the
restricted cutoffs. Any further zeros after restriction lie in $H|_E$,
so the tangential factors retain the exponents in \eqref{eq:mellin-factors}.

Choose smooth Hermitian metrics on the divisor line bundles of the
components of $B$. Replacing a coordinate weight $|z_G|^{2v_G}$ on $E$
by the norm of the canonical section to the power $2v_G$ multiplies the
coefficient by a smooth factor, entire in $v_G$ and equal to one at
$v_G=0$. The tangential distribution is holomorphic at the origin, so
these factors leave its value there unchanged. We may set $w=0$ in the
coefficient and evaluate along $v_G=N_G\lambda$ at $\lambda=0$, also
when disjoint components of $G|_E$ share a parameter.
Since $\chi\circ\mu=1$ on $E$ and the restricted partition sums to one,
the tangential values sum to
$c_{n-1}\operatorname{Reg}\!\int_E\rho\wedge\overline{h_E\rho}$
by Lemma~\ref{lem:pairing}. Multiplication by $\pi/N_E$ gives
\eqref{eq:mellin-residue}, whose right side is nonzero by
\eqref{eq:test-pairing}.
\end{proof}

\begin{lemma}[Bernstein--Sato root at a point center]\label{lem:point-root}
If a retained divisor $G$ has center a closed point $x$ and slope
$\alpha=\nu_G/N_G$, then $b_{f,x}(-\alpha)=0$ and $b_f(-\alpha)=0$.
\end{lemma}

\begin{proof}
Take $E$ to be the strict transform of $G$ in the construction above.
Apply the local Bernstein identity $b_{f,x}(s)f^s=Q(s)f^{s+1}$ on $U$
to the holomorphic exponent $s_1=s$, as in
\cite[proof of Theorem~5.3.2]{Ig}. For $\Re s\gg0$, integration by
parts gives
\begin{equation}\label{eq:mellin-bernstein}
 b_{f,x}(s) Z_h(s)
 =c_n\int_U f|f|^{2s}\overline h\,Q(s)^{\mathrm t}\chi\,
                  \omega\wedge\overline\omega ,
\end{equation}
where the transpose is taken with respect to $\omega$. To justify
integration across $D$, put $r=\operatorname{ord}Q$. For
$2\Re s>a+r$, the function $w\mapsto w|w|^{2s}/\bar w^a$ extends
by zero to a $C^r$ function on $\C$. Indeed, on the punctured disc
its real derivatives of order $j\le r$ are
$O_s(|w|^{2\Re s+1-a-j})$. Composing with $f$ and multiplying by
$\overline P$ gives a $C^r$ extension of $f|f|^{2s}\overline h$.
Off $D$ the holomorphic derivatives in $Q(s)$ annihilate the factor
$\bar f^{\,s}\overline h$ in local branches. The resulting identity
$Q(s)(f|f|^{2s}\overline h)=b_{f,x}(s)|f|^{2s}\overline h$ is
single-valued and extends by continuity across $D$. Transposing
against $\chi$ gives \eqref{eq:mellin-bernstein} on that half-plane.

The test function $Q(s)^{\mathrm t}\chi$ depends polynomially on $s$
and has support in $\Supp\chi$, so the charts in the proof of
Lemma~\ref{lem:zeta-residue} also cover the right side. Along each
boundary prime $F$ the factor $f$ increases the holomorphic exponent
in \eqref{eq:mellin-factors} to
\[
        \delta_F+N_F=\nu_F+(1-\alpha)N_F>0,
\]
since $X$ is smooth and $\alpha\le1$ by \eqref{eq:slope-bound}.
Thus no polar hyperplane in \eqref{eq:angular-poles} passes through
$\mathbf v=0$. Lemma~\ref{lem:local-mellin} gives meromorphic
continuation of the shifted integral along $v_F=N_F(s+\alpha)$
and holomorphy at $s=-\alpha$.
Taking residues in the continued identity \eqref{eq:mellin-bernstein}
gives $b_{f,x}(-\alpha)\Res_{s=-\alpha}Z_h(s)=0$ because the pole
of $Z_h$ is simple. Its residue is nonzero by Lemma~\ref{lem:zeta-residue},
so $b_{f,x}(-\alpha)=0$. Every algebraic Bernstein identity analytifies
at $x$, hence $b_{f,x}$ divides $b_f$.
\end{proof}

\section{Transverse slices and monodromy}\label{sec:slice}

We pass to transverse slices to prove Theorem~\ref{thm:retained} for arbitrary centers and deduce Theorem~\ref{thm:main}.

A local analytic form of Lemma~\ref{lem:slice-root} follows from Bahloul
\cite[Theorem~1(iii) and Remark~1(d)]{Bah}. For a family of germs
along a section over a smooth parameter germ, the local polynomial
$b_{\rm gen}(s)$ of a general fiber at the section divides
$\prod_{j=0}^N b_{\rm abs}(s+j)$ for some $N$, where $b_{\rm abs}$
is the polynomial of the total germ at the marked point of the
central fiber. Bahloul derives this divisibility from an iterated
functional equation of Biosca. Since Bernstein--Sato roots are strictly
negative, a root of $b_{\rm gen}$ in $[-1,0)$ must occur in
$b_{\rm abs}$. We give the argument for entire affine fibers over a
common Zariski-open set of parameters.

\begin{lemma}\label{lem:slice-root}
Let $X$ be smooth affine with an \'etale coordinate system
$(u_1,\dots,u_k,y_1,\dots,y_l)$. Let $X_a$ denote the fiber of the
map $u:X\to\mathbb A^k$ at $a$, and let $f_a=f|_{X_a}$. Suppose $f$
is nonzero on each component of a general fiber. Fix $0<\alpha\le1$.
If $b_f(-\alpha)\ne0$, there is a nonempty open
$V\subset\mathbb A^k$ such that $b_{f_a}(-\alpha)\ne0$ for every
$a\in V(\C)$.
\end{lemma}

\begin{proof}
Write $R=\OO_X(X)$, $K=\C(u_1,\ldots,u_k)$, $R_K=R\otimes_{\C[u]}K$, and
$\mathcal A=\C[s]_{(s+\alpha)}$. On the affine $X$, exactness of
global sections on quasi-coherent modules gives a global Bernstein
operator. Dividing by $b_f(s)$ gives
$f^s=P(s)f^{s+1}$ with $P\in D_{R/\C}\otimes_{\C}\mathcal A$ of some
order $r$, where $D_{R/\C}$ denotes the ring of differential operators.

We also need an identity over $R_K$ using only the relative
derivatives $\partial_y$. Fix an abstract isomorphism
$\overline K\simeq\C$, which exists because both fields are algebraically
closed of characteristic zero and have the same infinite transcendence
degree over $\Q$. Transporting scalars identifies algebraic differential
operators over these fields. On the resulting smooth complex affine
variety, the least common multiple of the local polynomials gives a
common Bernstein polynomial $b(s)\in\Q[s]$ with strictly negative
roots by Kashiwara \cite[Theorem~3.3 and Corollary~5.2]{Kas76}
(use $1$ on a component where $f$ is invertible). In the \'etale coordinates its
operator is a finite sum of monomials in $\partial_y$. The identity
is linear in their coefficients, so applying a $K$-linear retraction
$\overline K\to K$ fixing $1$ descends it to $R_K$.

For $1\le j\le r$ the polynomials $b(s+j)$ are units in $\mathcal A$,
since $j-\alpha\ge0$. Iterating the relative identity gives
$f^{s+1}=Q_r(s)f^{s+r+1}$ with
$Q_r\in D_{R_K/K}\otimes_{\C}\mathcal A$. Extend $P$ to $R_K$ by
localization. In $PQ_r$ commute the parameter derivatives to the
right. The resulting expression is $\sum_{|I|\le r}R_I(s)\partial_u^I$
with relative $R_I(s)$, since the coordinate derivations commute and
$\partial_u$ preserves $R_K$. For $|I|\le r$, the chain rule gives
$\partial_u^I f^{s+r+1}\in R_K[s]f^{s+1}$, since each derivative
lowers the exponent by at most one. Thus $f^s=PQ_r f^{s+r+1}$ lies
in $(D_{R_K/K}\otimes_{\C}\mathcal A)f^{s+1}$. Clearing the scalar
denominators gives
\[
      c(s)f^s=H(s)f^{s+1},\qquad c\in\C[s],\quad
      c(-\alpha)\ne0,\quad H\in D_{R_K/K}[s].
\]
Invert one nonzero $g\in\C[u]$ to put every coefficient of $H$ in
$R_g[s]$ in its $\partial_y$-expansion. Smoothness of $u$ gives
$R_g[1/f,s]\hookrightarrow R_K[1/f,s]$, so the identity already holds
over $R_g$. The relative derivatives and their chain-rule action on
$f^s$ specialize at each $a$ with $g(a)\ne0$. On shrinking to the open
in the hypothesis, $b_{f_a}$ consequently divides $c$.
\end{proof}

\begin{proof}[Completion of the proof of Theorem~\ref{thm:retained}]
The slope bound is \eqref{eq:slope-bound}, and
Lemma~\ref{lem:point-root} gives both root assertions for a point center. Suppose
$Z=\pi(G)$ has dimension $k>0$. Near a general smooth point of $Z$,
choose a smooth affine open $U\subset X$ with \'etale coordinates
$(u,y)$ for which $u|_{Z\cap U}$ is dominant and generically finite.
Restrict the models to $U$, write $W_a$ for the scheme-theoretic
fiber of any model $W$ over $a\in\mathbb A^k(\C)$, and take a log
resolution $p:Y\to S$ which resolves $(U,D_{\rm red}|_U)$ and is
an isomorphism off $D$.

We work on a single dense open of the base. First, $Z\cap U_a$ is
nonempty and finite there, and $D_a$ contains no component of $U_a$.
The generic fibers of $S$, $B_S$, and its primes are obtained by
localization, so the first is normal and the others reduced.
These properties are geometric over the perfect function field
$\C(u_1,\ldots,u_k)$. By
generic flatness and constructibility \cite[Proposition~8.9.4 and Corollary~9.9.5]{EGA},
we may take $S_a$ normal and all indicated boundary fibers reduced.
The flat fibers of $S$ and of each horizontal boundary prime are
cut out by a regular sequence of $k$ base parameters, so are pure of
dimensions $\dim U-k$ and $\dim U-k-1$, respectively.
On $Y$, generic smoothness applied to $Y$ and all intersections of
its boundary components makes the nonempty slices smooth of the
expected dimensions. The induced boundary on $Y_a$ therefore has
simple normal crossings. Vertical components can be omitted throughout.

On $S$ the singular loci of $S$ and its boundary primes, and the
intersections of distinct boundary primes, have codimension at
least two. By the fiber dimension theorem their slices have
codimension at least two in $S_a$ after further shrinking. Generic
smoothness on the smooth loci of $S$ and of each horizontal prime
then makes $S$ smooth over the base at all codimension-one points
of $S_a$, with at most one boundary prime there, smooth and Cartier
relative to the base. Thus $(B_S)_a$ is the reduced inverse image
of $D_a$, and $Y_a\to S_a\to U_a$ are projective and birational on
corresponding components, since they are isomorphisms off these
divisors and their inverse images.

Use the coordinate volume form on $U$ and its rational pullbacks.
On the relatively smooth loci of $S$ and $Y$, divide by
$du_1\wedge\cdots\wedge du_k$ using the cotangent determinant
sequence. The relative top forms restrict to rational canonical
forms on the fibers. Along the boundary in codimension one the
relative form is a unit times a power of a boundary equation times
a relative volume form. This follows since its divisor is supported
over $D$ and a single smooth prime passes through the generic point
of each sliced boundary component. Its equation restricts to a
uniformizer there. Consequently the orders of the volume form and
of $f$ along a horizontal boundary prime are preserved on each
component of its slice. For a Cartier multiple $m(K_S+B_S)$,
normality of $S_a$ now identifies its restriction with
$m(K_{S_a}+(B_S)_a)$ by equality in codimension one. This multiple
remains relatively ample. Writing $K_Y+\Gamma=p^*(K_S+B_S)$, the
same calculation on $Y_a$ and base change for Cartier pullback give
\[
       K_{Y_a}+\Gamma|_{Y_a}
       =p_a^*\bigl(K_{S_a}+(B_S)_a\bigr).
\]
On the log resolution $p_a$, the divisor $\Gamma|_{Y_a}$ has
simple normal crossings support. Its coefficients are inherited from
$\Gamma$ and hence are at most one, since $(S,B_S)$ is log canonical.
Thus $(S_a,(B_S)_a)$ is log canonical by the discrepancy criterion.
It follows that $S_a\to U_a$ is the log canonical modification on each component.

For every $a$ in this dense open, $G_a$ is nonempty and each of its
components $G'$ is a retained divisor centered at a point of $Z\cap U_a$,
with $(N_{G'},\nu_{G'})=(N_G,\nu_G)$. Lemma~\ref{lem:point-root},
applied on the smooth component of $U_a$ containing the center of $G'$,
gives the root $-\nu_G/N_G$ at that center. Since a Bernstein identity
on $U_a$ restricts to the germ there, $b_{f_a}(-\nu_G/N_G)=0$.
Lemma~\ref{lem:slice-root} forces $b_{f|_U}(-\nu_G/N_G)=0$,
since two dense opens of the base meet. A Bernstein identity on $X$
restricts to $U$, giving the same root of $b_f$.
\end{proof}

\begin{proof}[Proof of Theorem~\ref{thm:main}]
By Theorem~\ref{thm:witness}, $\alpha=\nu_G/N_G$ for a retained
divisor $G$, so Theorem~\ref{thm:retained} gives $b_f(-\alpha)=0$. The nearby cycle and Milnor monodromy assertions are consequences of a classical theorem of Malgrange and Kashiwara \cite[Theorem~3.4 and \S7, pp.~263--265]{Mal} and \cite[\S7, Theorem~2 and Remark~2]{Kas83}.%, every analytic neighborhood of $z$ contains a point $y\in D$ at which $\exp(-2\pi i\alpha)$ is an eigenvalue on $H^j(F_y,\C)$ for some $j$. For the nearby-cycle comparison in terms of the $V$-filtration, apply Kashiwara \cite[\S7, Theorem~2 and Remark~2]{Kas83} to the direct image of $\OO_X$ under the graph embedding of $f$, with respect to the coordinate $t$ on $X\times\mathbb A^1$.
\end{proof}

\end{document}